\documentclass[11pt]{amsart}

\usepackage{amssymb}
\usepackage[all]{xy}
\begin{document}
\title{Counting conics in complete intersections}
\author{Laurent BONAVERO, Andreas H\"ORING}
\date{September 22, 2009}
\maketitle
\noindent
\def\restriction{\string |}
\newcommand{\pp}{\rm ppcm}
\newcommand{\pg}{\rm pgcd}
\newcommand{\Ker}{\rm Ker}
\newcommand{\C}{{\mathbb C}}
\newcommand{\Q}{{\mathbb Q}}
\newcommand{\GL}{\rm GL}
\newcommand{\SL}{\rm SL}
\newcommand{\diag}{\rm diag}
\newcommand{\N}{{\mathbb N}}
\def\finpreuve
{\hskip 3pt \vrule height6pt width6pt depth 0pt}

\newtheorem{lemma1}{}[section]

\newenvironment{lemm}{\begin{lemma1}{\bf Lemma.}}{\end{lemma1}}
\newenvironment{example}{\begin{lemma1}{\bf Example.}\rm}{\end{lemma1}}
\newenvironment{question}{\begin{lemma1}{\bf Question.}\rm}{\end{lemma1}}
\newenvironment{abs}{\begin{lemma1}\rm}{\end{lemma1}}
\newenvironment{theo}{\begin{lemma1}{\bf Theorem.}}{\end{lemma1}}
\newenvironment{prop}{\begin{lemma1}{\bf Proposition.}}{\end{lemma1}}
\newenvironment{cor}{\begin{lemma1}{\bf Corollary.}}{\end{lemma1}}
\newenvironment{rem}{\begin{lemma1}{\bf Remark.}\rm}{\end{lemma1}}
\newenvironment{defi}{\begin{lemma1}{\bf Definition.}}{\end{lemma1}}
\newenvironment{construction}{\begin{lemma1}{\bf Construction.}}{\end{lemma1}}
\newenvironment{conjecture}{\begin {lemma1}{\bf Conjecture.}}{\end{lemma1}}
\newenvironment{problem}{\begin{lemma1}{\bf Problem.}}{\end{lemma1}}

\newcommand{\CC}{{\mathbb C}}
\newcommand{\ZZ}{{\mathbb Z}}
\newcommand{\RR}{{\mathbb R}}
\newcommand{\QQ}{{\mathbb Q}}
\newcommand{\FF}{{\mathbb F}}
\newcommand{\PP}{{\mathbb P}}
\newcommand{\NN}{{\mathbb N}}
\newcommand{\codim}{\operatorname{codim}}
\newcommand{\Ho}{\operatorname{Hom}}
\newcommand{\Pic}{\operatorname{Pic}}
\newcommand{\NE}{\operatorname{NE}}
\newcommand{\Nun}{\operatorname{N}}
\newcommand{\card}{\operatorname{card}}
\newcommand{\Hilb}{\operatorname{Hilb}}
\newcommand{\mult}{\operatorname{mult}}
\newcommand{\vol}{\operatorname{vol}}
\newcommand\sO{{\mathcal O}}
\newcommand{\divi}{\operatorname{div}}
\newcommand{\pr}{\operatorname{pr}}
\newcommand{\con}{\operatorname{cont}}
\newcommand{\ima}{\operatorname{Im}}
\newcommand{\rk}{\operatorname{rk}}
\newcommand{\Exc}{\operatorname{Exc}}
\newcounter{subsub}[subsection]
\def\thesubsub{\thesubsection .\arabic{subsub}}
\def\subsub#1{\addtocounter{subsub}{1}\par\vspace{3mm}
\noindent{\bf \thesubsub ~ #1 }\par\vspace{2mm}}
\def\coker{\mathop{\rm coker}\nolimits}
\def\pr{\mathop{\rm pr}\nolimits}
\def\im{\mathop{\rm Im}\nolimits}
\def\hfl#1#2{\smash{\mathop{\hbox to 12mm{\rightarrowfill}}
\limits^{\scriptstyle#1}_{\scriptstyle#2}}}
\def\vfl#1#2{\llap{$\scriptstyle #1$}\big\downarrow
\big\uparrow
\rlap{$\scriptstyle #2$}}
\def\diagram#1{\def\normalbaselines{\baselineskip=0pt
\lineskip=10pt\lineskiplimit=1pt}   \matrix{#1}}
\def\limind{\mathop{\oalign{lim\cr
\hidewidth$\longrightarrow$\hidewidth\cr}}}

\long\def\InsertFig#1 #2 #3 #4\EndFig{
\hbox{\hskip #1 mm$\vbox to #2 mm{\vfil\includegraphics{#3}}#4$}}
\long\def\LabelTeX#1 #2 #3\ELTX{\rlap{\kern#1mm\raise#2mm\hbox{#3}}}


\newcommand{\merom}[3]{\ensuremath{#1:#2 \dashrightarrow #3}}
\newcommand{\holom}[3]{\ensuremath{#1:#2  \rightarrow #3}}
\newcommand{\fibre}[2]{\ensuremath{#1^{-1} (#2)}}
\newcommand{\Z}{{\mathbb Z}}
\newcommand{\chow}[1]{\ensuremath{\mathcal{C}(#1)}}
\newcommand{\nons}{\operatorname{nons}}
\newcommand\sI{{\mathcal I}}

\setcounter{tocdepth}{1}


{\let\thefootnote\relax
\footnote{
\textbf{A.M.S.~classification :} 14N10, 14M10. 
}}

\vspace{-1cm}


\begin{center}
\begin{minipage}{130mm}
\scriptsize

{\bf Abstract.} We count the number of conics through two general
points in complete intersections when this number is finite
and give an application in terms of quasi-lines.
\end{minipage}
\end{center}

\section{Introduction}

Let $X$ be a complex projective manifold of dimension $n$.
A quasi-line $l$ in $X$ is a {\it smooth} rational curve $f~: \PP^1 
\hookrightarrow X$
such that $f^*T_X$ is the same as for a line in $\PP^n$, {\em i.e.} 
is isomorphic to 
$${\mathcal O}_{\PP^1}(2) \oplus {\mathcal O}_{\PP^1}(1)^{\oplus n-1}.$$

Let $X$ be a smooth projective variety containing a quasi-line $l$. 
Following Ionescu and Voica \cite{IV03}, we 
denote by $e(X,l)$ the number of quasi-lines which are deformations
of $l$ and pass through two given general points of $X$. We 
denote by $e_0(X,l)$ the number of quasi-lines which are deformations
of $l$ and pass through a general point $x$ of $X$ with
a given general tangent direction at $x$. Note that one always has 
$e_0(X,l) \leq e(X,l)$, but in general the inequality may be strict \cite[p.1066]{IN03}.

\begin{theo} \label{maintheorem} 
Let $X \subset \mathbb P ^{n+r}$ be a {\em general} smooth $n$-dimensional
complete intersection of multi-degree $(d_1,\ldots,d_r)$. 
Assume moreover that 
$$ d_1+\cdots+d_r = \frac{n+1}{2} +r.$$
Then
\begin{enumerate}

\item the family of conics contained in $X$ is a nonempty, 
smooth and irreducible
component of the Chow scheme $\chow{X}$,
\item a general conic $C$ contained in $X$ is a quasi-line
of $X$ and
$$
\displaystyle{
e_0(X,C)=e(X,C)= \frac{1}{2}\prod _{i=1}^r (d_i -1)! d_i !}.
$$
\end{enumerate}
\end{theo}
 
The numerical assumption $d_1+\cdots+d_r = (n+1)/2 +r$ assures that
if $C$ is a conic in $X$, then 
$-K_X \cdot C = n+1$. This numerical condition is of course necessary
for a curve to be a quasi-line. Note that varieties appearing in 
our theorem are Fano varieties of dimension $n$ and index $(n+1)/2$; they
are well known to be the boundary Fano varieties with Picard number one
being conic-connected (see \cite{IR07}, Theorem 2.2). 

Using a degeneration argument, one can strengthen parts of the statement.

\begin{cor} 
Let $X \subset \mathbb P ^{n+r}$ be a smooth $n$-dimensional
complete intersection of multi-degree $(d_1,\ldots,d_r)$. 
If $ d_1+\cdots+d_r = (n+1)/2 +r$, the variety
$X$ contains a conic that is a quasi-line.
\end{cor}

By a theorem of Ionescu \cite{Ion05}, we obtain an immediate 
application of the theorem to formal geometry. Before stating it, let
us recall that a subvariety $Y$ of a variety $X$ is ${\rm G3}$ in $X$
if the ring $K(X_{|Y})$ of formal-rational functions of
$X$ along $Y$ is equal to $K(X)$. 

\begin{cor} \label{corollaryformal} 
Let $X \subset \mathbb P ^{n+r}$ be a general smooth $n$-dimensional
complete intersection of multi-degree $(d_1,\ldots,d_r)$
such that $d_1+\cdots+d_r = (n+1)/2 +r.$ Then any general
conic $C$ contained  in $X$ is ${\rm G3}$ in $X$.
In particular, if $(X,C)$ and $(X',C')$ are two such pairs
such that the formal completions $X_{|C}$ and $X'_{|C'}$ 
are isomorphic as formal schemes, there exists an
isomorphism from $X$ to $X'$ sending $C$ to $C'$.
\end{cor}

When this note was almost finished, we learned from L. Manivel that
A. Beauville  had obtained the formula $\displaystyle{e(X,l)= 
\frac{1}{2}\prod _{i=1}^r (d_i -1)! d_i !}$
as a consequence of his computation of the
quantum cohomology algebra $H^*(X,\mathbb Q)$  of 
a complete intersection \cite{Bea95}. 
We provide here a completely elementary proof. We end this note 
by mentioning a similar question where no elementary proof seems
to be known. 
 
\medskip

{\em We want to thank Fr\'ed\'eric Han, Paltin Ionescu and Laurent Manivel
for their interest in our work and helpful discussions.
The first author express his warm gratitude to
the organisers of the Hano\"{\i} conference, with a 
special mention to H\'el\`ene Esnault, for their kind invitation.
We also thank the referee for numerous detailed remarks.}

\section{Proofs}

We start by explaining the enumerative argument in the simplest case. 

\subsection{A well known example} 

Suppose that $X=\{s=0\}$ is a smooth cubic threefold in $\PP^4$. 
A general conic $C$ in $X$ is a quasi-line \cite[Thm.3.2]{BBI00}.
The basic idea of our proof is that counting conics in $X$ 
through $p$ and $q$ can be reduced to counting $2$-planes $\pi$ 
through $p$ and $q$  
such that the restriction $s|_{\pi}$ is a product 
of a polynomial of degree two and some residual polynomial. 
We will explain how to do this in general below, in
the case of the cubic threefold we can use a geometric construction.

It is a classical fact that
the lines in $X$ form an irreducible smooth family of dimension two and 
that there are exactly six lines passing through a 
general point of $X$ \cite[Prop.1.7]{AK77}.
Fix now two general points $p$ and $q$ in $X$, then the line 
$[pq]$ intersects $X$ in a third point $u$. 
For every line $l \subset X$ through $u$ 
there exists a unique plane $\pi_l$ containing $l$ and $[pq]$. 
The intersection $X  \cap \pi_l$ is the union of $l$ and a residual conic $C$.
Since $l$ does not pass through $p$ and $q$, the conic $C$
passes through $p$ and $q$. {\em Vice versa} 
the linear span of a conic $C \subset X$ passing through $p$ and $q$ is a 
$2$-plane $\pi_C$ containing the line $[pq]$.  
Since $C$ does not pass through $u$, the residual line passes through $u$.
Thus the conics through $p$ and $q$ are in bijection 
with the lines through $u$, so $e(X, C)=6$.

\medskip

Suppose now that we are in the general 
situation of Theorem \ref{maintheorem}.
We always assume that 
$X \subset \PP^{n+r}$ is a general 
smooth $n$-dimensional
complete intersection of multi-degree $(d_1,\ldots,d_r)$ 
with $d_i \geq 2$ for all $i$ and
$$ 
d_1+\cdots+d_r = \frac{n+1}{2} +r.
$$
Let $l \subset X$ be a smooth rational curve contained in 
$X$. Then 
$$ -K_X \cdot l = (n+r+1- (d_1+\cdots+d_r))\deg (l)= \frac{n+1}{2}\deg (l)$$
therefore $-K_X \cdot l = n+1$ if and only if $l$ is a conic.

\subsection{The main step}

{\em For any general points 
$p$ and $q$ of $X$, there exists
a conic contained in $X$ passing through $p$ and $q$.}

\medskip

Fix two distinct points in $\PP^{n+r}$, say
$p=[1:0:\cdots:0]$ and $q=[0:0:\cdots:1]$.
Suppose that $X$ is a general complete
intersection with equations
$$(s_1=0)\cap (s_2=0)\cap \cdots \cap (s_r=0)$$
passing through $p$ and $q$, 
where each $s_i \in H^0(\PP^{n+r},\sO_{\PP^{n+r}}(d_i))$
is general among sections vanishing at $p$ and $q$.

\medskip

Suppose there is a conic $C$ contained in $X$, passing
through $p$ and $q$ and let
$\pi _C$ the projective $2$-plane generated by $C$. 
If $s_C$ denotes the equation
defining $C$ in $\pi_C$, there exists for each $i=1,\ldots,r$ a 
$\tilde s _i \in H^0(\pi _C,\sO_{\pi _C}(d_i -2))$ (defining the residual
curve)
such that 
$$ 
(s_i)_{|\pi _C} = s_C \cdot \tilde s _i.
$$
Since $X$ is general, it does not contain the $2$-plane $\pi _C$
\cite[Thm. 2.1]{DM98}. Therefore $(s_i)_{|\pi _C}$
and $\tilde s _i$ are not zero for at least one $i$.

Conversely, let $\pi$ be a projective $2$-plane containing $p$ and $q$
and assume there exists a non-zero $s_C \in H^0(\pi,\sO_{\pi}(2))$
vanishing at $p$ and $q$ 
and, for each $i=1,\ldots,r$, there exists a 
$\tilde s _i \in H^0(\pi,\sO_{\pi}(d_i -2))$
such that 
$$ (s_i)_{|\pi} = s_C \cdot \tilde s _i,$$
then the conic $(s_C =0)$ is obviously contained in $X$.

Consider now the projective space of dimension $n+r-2$ parametrizing the
projective $2$-planes
in $\PP^{n+r}$ containing $p$ and $q$. 
Fixing homogeneous coordinates $[a_1 : \cdots : a_{n+r-1}]$ on this space,
let 
$$
\pi_{[a_1:\cdots:a_{n+r-1}]} 
= \{[x:za_1:\cdots:za_{n+r-1}:y]\mid [x:z:y]\in \mathbb P^2 \}
$$
be such a $2$-plane. Then 
$$ 
(s_i)_{|\pi_{[a_1:\ldots:a_{n+r-1}]}}(x,z,y)= 
\sum_{k=0}^{d_i} \sum_{a=0}^k s^i_{a,k}x^ay^{k-a}z^{d_i-k}
$$
where $s^i_{a,k}$ is a homogeneous polynomial
of degree $d_i-k$ in the variables $a_1,\ldots,a_{n+r-1}$.
The equation of an irreducible conic in this plane that passes through $p$ and $q$ is
$$
s_C = s_2z^2+s_1xz+s' _1yz+ xy.
$$
So for each $i=1,\ldots,r$, 
the equation $(s_i)_{|\pi} = s_C \cdot \tilde s _i$
can be written explicitly  
$$
\sum_{k=0}^{d_i} \sum_{a=0}^k s^i_{a,k}x^a y^{k-a} z^{d_i-k}
= (s_2z^2+s_1xz+s' _1yz+ xy)\times 
\sum_{k=0}^{d_i-2}\sum_{a=0}^k \tilde{s}^i_{a,k}x^ay^{k-a}z^{d_i-2-k}.
$$
Thus we have to solve the equations 
$$
s^i_{a,k} = s_2 \tilde{s}^i_{a,k}+s_1\tilde{s}^i_{a-1,k-1}
+s'_1\tilde{s}^i_{a,k-1}+\tilde{s}^i_{a-1,k-2}
$$
for any $0\leq k \leq d_i$ and $0\leq a \leq k$.

\medskip

Let us first solve this system (whose unknown variables are
$s_2$, $s_1$, $s'_1$ defining the conic and the $\tilde{s}^i_{a,k}$'s defining the residual curve) for each $i$ separately. 
Note that
$X$ passes through $p$ and $q$ if and only if $s^i_{0,d_i}= s^i_{d_i,d_i}=0$.
Therefore writing the $d_i-1$ equations $s^i_{a,d_i}=\tilde{s}^i_{a-1,d_i-2}$
for $1\leq a \leq d_i-1$ provides the $\tilde{s}^i_{a-1,d_i-2}$'s.

Considering the equations corresponding to $(a,k) = (0,d_i-1)$
and $(d_i-1,d_i-1)$ 
allows to find $s_2$ and $s_1$.  
Considering then the equations corresponding to $(a,k) = (1,d_i-1)$
and $(0,d_i-2)$ gives $\tilde{s}^i_{0,d_i-3}$
and $s'_1$ (in particular this determines the conic, 
if it exists~!).
Write down successively the equations 
for $(a,k)$, $a=1,\ldots,k-1$, $k=d_i-1,\ldots,2$ to find
all the $\tilde{s}^i_{a,k}$'s (this determines the residual curve 
$(\tilde{s}^i=0)$ !).

Therefore, the $r$ systems have a common solution if 
and only if the remaining
equations for each system are satisfied and the 
corresponding conic is the same for each $i$.
For each $i$, the remaining equations are ``universal formulas''
(meaning the coefficients just depend on the equations defining $X$) 
corresponding to $(a,k)=(0,d_i-3),\ldots,(0,0)$
and $(a,k)=(d_i-2,d_i-2),\ldots,(1,1)$. This gives
$2d_i-4$ equations of respective degrees $3,\ldots,d_i$ and 
$2,3,\ldots,d_i-1$ in the variables $a_1,\ldots,a_{n+r-1}$. 
The $3r-3$ equations saying that the conic is the
same for each $i=1,\ldots,r$ are $2r-2$ equations of degree $1$ and $r-1$
equations of degree $2$ in the variables $a_1,\ldots,a_{n+r-1}$.

Altogether, using the relation $d_1+\cdots+d_r = (n+1)/2 +r$,
this gives exactly $n+r-2$ homogeneous equations 
in the variables $a_1,\cdots,a_{n+r-1}$. We therefore get at least 
one solution. Moreover since $X$ is general, the coefficients 
$s^i_{a,k}$ appearing in the initial equations are general. Since they completely determine
the remaining $n+r-2$ homogeneous equations, these equations are general.
Thus the space of solutions is smooth and of the expected dimension, so
there are exactly $\displaystyle{
\frac{1}{2}\prod _{i=1}^r (d_i -1)! d_i !}$ solutions by Bezout's theorem.

\medskip

Let us briefly indicate how the same method gives 
{\em the number of conics contained in $X$, passing through 
$p$ and tangent to the line $(pq)$}. With the above notations,
we have $s^i_{d_i-1,d_i}= s^i_{d_i,d_i}=0$
and
we have to solve the $r$ systems
$$\sum_{k=0}^{d_i} \sum_{a=0}^k s^i_{a,k}x^a y^{k-a} z^{d_i-k}
= (s_2z^2+s_1xz+s' _1yz+ y^2)\times 
\sum_{k=0}^{d_i-2}\sum_{a=0}^k \tilde{s}^i_{a,k}x^ay^{k-a}z^{d_i-2-k}$$
which means 
$$s^i_{a,k} = s_2 \tilde{s}^i_{a,k}+s_1\tilde{s}^i_{a-1,k-1}
+s'_1\tilde{s}^i_{a,k-1}+\tilde{s}^i_{a,k-2}$$
for any $0\leq k \leq d_i$ and $0\leq a \leq k$. The remaining details are 
left to the reader.

\subsection{The space of conics is irreducible}
Let $\mathbb{G}(2,n+r)$ be the Grassmannian of projective $2$-planes contained
in $\PP^{n+r}$ and $E$ be the tautological rank $3$-bundle 
on $\mathbb{G}(2,n+r)$.
The Hilbert scheme of conics in $\PP^{n+r}$ is 
the projectivisation\footnote{
In this article we follow the convention that the projectivisation of a vector bundle $E$ is the variety of lines of $E$.}
of $S^2 E^*$. Denote by $\varphi: \mathbb P (S^2 E^* ) \to \mathbb G(2,n+r)$
the natural map. We have an exact sequence on $\mathbb P (S^2 E^* )$~:
$$ 
(*) \qquad 0 \to \bigoplus_{i=1}^r \varphi^*S^{d_i-2} E^*  \otimes 
\mathcal O_{\mathbb P (S^2 E^* )}(-1) 
\to  \bigoplus_{i=1}^r \varphi^*S^{d_i} E^*
\to \mathcal Q \to 0 
$$
defining a vector bundle $\mathcal Q$ of rank $n+1+3r$.
Since $X$ is a complete
intersection $(s_1=0)\cap (s_2=0)\cap \cdots \cap (s_r=0)$,
the $s_i$'s induce by restriction to $2$-planes,
pull-back and projection onto $\mathcal Q$ a section 
of $\mathcal Q$ whose zero locus $Z$ is precisely the set of
conics contained in $X$. Since $E^*$ is globally generated, 
the images of sections $(s_1, \ldots, s_r)$ 
give a vector space $V \subseteq H^0(\mathbb P (S^2 E^* ), \mathcal Q)$
that globally generates $\mathcal Q$.
Applying Bertini's theorem to this subspace we see that 
the zero locus of a general section in $V$ is smooth. 
Since $X$ is supposed to be a {\em general} complete intersection,
$Z$ is smooth and proving its irreducibility reduces to showing
that $h^0(Z,\mathcal O _Z)=1$. 
By the Koszul resolution of $ \mathcal O_Z$, 
it is enough to show that for any $1\leq j \leq \rk \mathcal Q$ 
$$
h^j(\mathbb P (S^2 E^* ), \wedge^j \mathcal Q^*)=0.
$$
Using the exact sequence $(*)$, this easily reduces to showing 
that 
for any $1\leq j \leq \rk \mathcal Q$ and any $0\leq k \leq j$,
$$ 
H^k(\mathbb P (S^2 E^* ),\wedge^k (\oplus_{i=1}^r \varphi^*S^{d_i} E)
\otimes S^{j-k}( (\oplus_{i=1}^r \varphi^*S^{d_i-2} E)  \otimes 
\mathcal O_{\mathbb P (S^2 E^* )}(1)))=0.
$$
Since the higher direct images with respect to $\varphi$ vanish, 
it is sufficient
to show that for any 
$1\leq  j \leq \rk \mathcal Q$ and for any $0\leq k \leq j$,
we have
$$ 
 H^k(\mathbb G (2,n+r),\wedge^k (\oplus_{i=1}^r S^{d_i} E)
\otimes S^{j-k}( (\oplus_{i=1}^r S^{d_i-2} E)\otimes S^{2} E )))=0.
$$

This will follow from Bott's theorem
applied on $\mathbb{G}(2,n+r)$. Indeed, using Schur functor notation,
let $S_b E$ be an irreducible factor appearing in the decomposition
of $\wedge^k (\oplus_{i=1}^r S^{d_i} E)
\otimes S^{j-k}( (\oplus_{i=1}^r S^{d_i-2} E)\otimes S^{2} E )$
where $b=(b_1,b_2,b_3)$ is a triple of integers $b_1\geq b_2 \geq b_3 \geq 0$. 
By the  Littlewood-Richardson rule,
we get 
$$
(**) \qquad
b_2 + b_3 \geq k-r \mbox{ and } b_3 \geq k-(d_1+\ldots+d_r)-r
= k-(n+1)/2-2r.
$$

On the other hand by Bott's theorem,
the whole cohomology of $S_b E$ vanishes except
maybe in the following cases~: 
\begin{enumerate} 
\item $k=n+r-2$ and $(b_1,b_2,b_3)=(b_1,0,0)$ with $b_1 \geq n+r-1$,
\item $k=n+r-2$ and $(b_1,b_2,b_3)=(b_1,1,0)$ with $b_1 \geq n+r-1$,
\item $k=n+r-2$ and $(b_1,b_2,b_3)=(b_1,1,1)$ with $b_1 \geq n+r-1$,
\item $k=2(n+r-2)$ with $b_2 \geq n+r$
and $b_3=0,1,2$,
\item $k=3(n+r-2)$ with $b_3 \geq n+r+1$.
\end{enumerate} 

The case  $n=3$ has been dealt with by B{\u{a}}descu, Beltrametti
and Ionescu \cite{BBI00}, so we may assume $n \geq 5$ since $n$ is odd.

In the first three cases, we
get $k-r =n-2 \geq 3 > b_2+b_3 = 0,1,2$, which is excluded by $(**)$.
In case (4), since $n\geq 5$, we get $k-(n+1)/2-2r = 3(n-3)/2 > b_3=0,1,2$, 
which is again excluded by $(**)$.
Case (5) is also excluded since  
we are only interested in the situation where
$k \leq \rk \mathcal Q = n+1+3r$, but $3(n+r-2)> n+1+3r$ when $n\geq 5$.

\medskip

We obtain the following corollary of the proof.

\begin{cor} 
Let $X \subset \mathbb P ^{n+r}$ be a general smooth $n$-dimensional
complete intersection of multi-degree $(d_1,\ldots,d_r)$. 
Assume moreover that 
$$ 
d_1+\cdots+d_r \leq \frac{n+1}{2} +r
$$
and $n \geq 5$.
Then the family of conics contained in $X$ is a nonempty, 
smooth and irreducible
component of the Chow scheme $\chow{X}$.
\end{cor}

Let us also mention that Harris, Roth and Starr have shown the irreducibility of
the space of smooth rational curves of arbitrary degree $e$ for general hypersurfaces of low degree $d$ \cite{HRS04}.

\subsection{Conics are quasi-lines}
By the first step, there exists a conic $C$ 
passing through two general points. Such a conic is necessarily smooth:
a line $d$ contained in $X$ and passing through a general point satisfies
\[
T_X|_d \simeq {\mathcal O}_{\PP^1}(2) 
\oplus {\mathcal O}_{\PP^1}(1)^{\oplus \frac{n-3}{2}}
\oplus {\mathcal O}_{\PP^1}^{\oplus \frac{n+1}{2}},
\]
so an easy dimension count shows that two general points are not 
connected by a chain of two lines.
Thus $C$ smooth and its deformations with a fixed
point cover a dense open subset in $X$. This implies that
the normal bundle $N_{C/X}$ is ample \cite[Prop.4.10]{Deb01} and since $-K_X\cdot C =n+1$,
the curve $C$ is a quasi-line. 

\subsection{Proof of the Corollary \ref{corollaryformal}}
The irreducibility of the variety of conics gives
$$
\displaystyle{
e_0(X,l)=e(X,l)= \frac{1}{2}\prod _{i=1}^r (d_i -1)! d_i !}.
$$
The equality $e_0(X,l)=e(X,l)$ implies that general conics 
are {\rm G3} in $X$ \cite[Cor. 4.6]{Ion05}, in particular 
\cite[Cor. 4.7, Cor.1.9]{Ion05} apply. 

\section{A similar question}

Using exactly the same method as developed in \S 2.3 ,
one can prove the following result, left to the reader. 

\begin{prop} 
Let $X_d \subset \mathbb P ^{n+1}$ be a general smooth $n$-dimensional
hypersurface of degree $d$. 
Then, for $n\ge 7$ and $d \leq n+1$, 
the family of conics contained in $X_d$ is a nonempty, 
smooth and irreducible
component of dimension $3n-2d+1$ of the Chow scheme $\chow{X_d}$.
\end{prop}

In the case of $d=n+1$, there is a finite number of conics passing
through a general point of $X_{n+1}$. Let us denote by $N_{n+1}$
this number. It seems that there are no known elementary 
method to compute this number. A general formula comes from
the calculation of some Gromov-Witten invariants using mirror
symmetry and an ordinary differential 
equation introduced by Givental. The following lines 
were written while reading \cite{JNS04} and \cite{Jin05}.

\begin{prop} \label{coniques} 
{\bf (Coates, Givental - Jinzenji, Nakamura, Suzuki)}
Let 
$X_n \subset \mathbb P ^n$ be a general smooth
hypersurface of degree $n$ in $\mathbb P ^n$. Let
$N_n$ be the number of conics passing
through a general point of $X_{n}$. 
Then
$$ N_n = \frac{(2n)!}{2^{n+1}}-\frac{(n!)^2}{2}.$$ 
\end{prop}

Let us briefly explain where this result comes from. If $a$, $b$,
$c$ et $d$ are four integers, let 
$\langle \mathcal O_a \mathcal O_b \mathcal O_c\rangle_d$
be the Gromov-Witten
invariant counting the number
(possibly infinite) of rational curves of degree $d$
contained in     
$X_n$ and meeting $3$ general subspaces of $\mathbb P ^n$, 
of respective codimension 
$a$, $b$ and
$c$. When $a$, $b$ or $c$ are equal to $1$, each such rational curve
has to be counted $d$ times since the intersection 
of a degree $d$ curve intersects a general hyperplane in $d$ points.
Since a general line meets $X_n$ in $n$ points, 
we get that 
$N_n = \langle \mathcal O _1 \mathcal O _1 \mathcal O _{n-1}\rangle _2 /4n$. 
In \cite{Jin05} are introduced some constants $\tilde L _m ^{n+1,n,d}$,
called ``structure constants of the quantum cohomology
ring of $X_n$''. They satisfy the following formula:
$$ \sum_{m=0}^{n-1}\tilde L _m ^{n+1,n,1}w^m=n\prod_{j=1}^{n-1}(jw+(n-j))$$
and 
$$ \sum_{m=0}^{n-2}\tilde L _m ^{n+1,n,2}w^m=
\sum_{j_2=0}^{n-2}\sum_{j_1=0}^{j_2}\sum_{j_0=0}^{j_1}
\tilde L _{j_1} ^{n+1,n,1}\tilde L _{j_2+1} ^{n+1,n,1}w^{j_1-j_0}\left(
\frac{1+w}{2}\right) ^{j_2-j_1}.$$
 
It is also shown in \cite{Jin05} that for every integer $m$, 
$0 \leq m \leq n-2$,
we have $\tilde L _m ^{n+1,n,2}= 
\langle \mathcal O _1 \mathcal O _{n-1-m} \mathcal O _{m+1}\rangle _2 /n$. 

Then the proposition follows by evaluating the $w^{n-2}$ 
coefficient in the second formula above,
the $w^{n-1}$ coefficient in the first
and putting $w=2$.


\medskip

{\em \small
\noindent Laurent Bonavero. 
Institut Fourier, UMR 5582, 
Universit\'e de Grenoble 1, BP 74. 
38402 Saint Martin d'H\`eres. France. 
\\
e-mail : laurent.bonavero@ujf-grenoble.fr

\noindent Andreas H\"oring. 
Universit\'e Paris 6, Institut de Math\'ematiques de Jussieu, Equipe de Topologie et G\'eom\'etrie Alg\'ebrique, 175, rue du Chevaleret, 75013 Paris, France.
\\
\noindent e-mail : hoering@math.jussieu.fr
}

\end{document}